\documentclass[]{gENO2e}

\begin{document}

\doi{10.1080/0305215X.YYYY.CATSid}
 \issn{1029-0273}
\issnp{0305-215X}
\jvol{00} \jnum{00} \jyear{2011} \jmonth{October}

\markboth{Xiaojun Zhou and Chunhua Yang}{canonical dual  algorithms}

% \articletype{GUIDE}

\title{Four strategies to develop canonical dual  algorithms for\\
 global optimization problems}

\author{Xiaojun Zhou$^{\dag}$$^{\ddag}$\thanks{$\dag$Corresponding author. Email: tiezhongyu2010@gmail.com
\vspace{6pt}} and Chunhua Yang$^{\ddag}$\\\vspace{6pt}  $^{\dag}${\em{School of Science, Information Technology and Engineering, University of Ballarat, Victoria 3353, Australia}};
$^{\ddag}${\em{School of Information Science and Engineering, Central South University, Changsha 410083, China}}\\\vspace{6pt}}
%\received{v3.9 released December 2010}

\maketitle

\begin{abstract}
The canonical duality theory has provided with a unified analytic solution to a range of discrete and continuous problems in global optimization, which can transform a nonconvex primal problem to a concave maximization dual  problem over a convex domain without duality gap. This paper shows that under certain conditions,  this canonical dual problem is equivalent to the standard semi-definite programming (SDP) problem, which can be solved by well-developed software packages. In order to avoid certain difficulties of using the SDP method,
 four strategies are proposed based on unconstrained approaches, which can be used to develop algorithms for solving some challenging problems. Applications are illustrated by  fourth-order polynomials benchmark optimization problems.\bigskip

\bigskip
\end{abstract}
\vspace{-42pt}
\bigskip\bigskip\bigskip
\section{Introduction and motivation}\label{intro}
Numerical  optimization methods  are usually categorized into deterministic and stochastic, both of them have found extensive  applications in real-world problems \citep{hendrix,shmoys}.
The deterministic methods such as Newton's method were considered as ``local search" because they are dependent on
initial point to a large extent and finally arrive at the ``neighborhood" of the initial point, while the stochastic
methods such as genetic algorithm were regarded as ``global search" due to their search ability on the whole space.
However, the premature convergence and easily getting trapped into local optimum are common phenomena for stochastic
 algorithms \citep{back,bentley}.
  On the other hand, some deterministic global optimization algorithms have been proposed
  in recent years, and they are able to solve much general optimization problems such as nonconvex continuous, mixed-integer,
   differential-algebraic,
  bilevel, and non-factorable problems \citep{floudas}.
  The development of deterministic and stochastic just indicate what the ``No Free Lunch Theorems" means,
  that is, there exists no algorithm which is better than its competitor over all problems \citep{wolpert}. \\

\indent In the meanwhile, a novel global optimization theory called the {\em canonical duality theory}  has been developed during the past 20 years. The kernels of the theory consist of a canonical dual transformation methodology, a complementary-dual principle, and a triality theory \citep{gao1,gao2}.  The main merit is that this theory can transform nonconvex/nonsmoonth/discrete optimization/variational problems into continuous concave maximization  problems over convex domains, which can be solved easily, under certain conditions,  by many well-developed algorithms and softwares. Therefore, the canonical duality theory has been used successfully for solving a large class of challenging problems in computational biology \citep{zgy}, engineering mechanics \citep{Gao-Sherali-AMMA09,gao-yu2008,santos-gao}, information theory \citep{latorre-gao}, network communications \citep{gao3}, nonlinear dynamical systems \citep{ruan-gao-ima},
and some NP-hard problems in global optimization \citep{fang-gaoetal07,gao2010,gao-watsonetal,wangetal}.  \\

\indent  However, it was realized that the canonical dual problem may have  no  critical point in the dual feasible space and in this case, the primal problem could be NP-hard \citep{gao-jimo07}.
By introducing a linear perturbation term to the primal problem or a quadratic perturbation term to the dual problem, the issue can be partially tackled with to some extent but is still an open problem \citep{wangetal}.
For one thing, it is not easy to find such an appropriate perturbation. For another, only approximate solution is obtained due to the perturbation. On the other hand, it is undoubtedly that solving a constrained optimization problem (a continuous concave maximization  problem over convex domain) is much more difficult than an unconstrained one. Therefore, some approaches, such as the penalty function method, aim to convert a constrained minimization problem into an equivalent unconstrained one to reduce the computational complexity.\\

\indent As is known to us, for nonconvex optimization problem, methods based on gradient are dependent on initial point, and choosing a good initial point can reach a good solution in the end. To overcome local optimality, it usually requires some type of diversification to find the global optimum. For instance, the multi-start methods, which are applied by starting from multiple random initial solutions, are widely used to realize diversification \citep{marti}. However, it remains hard to construct good initial solutions.\\

\indent Fortunately, the global optimality condition contained in triality theory, can identify the global minimum, which provides with greatly useful information to select a good initial point and can be utilized to develop related canonical dual algorithms. In this study, we show that, under certain conditions, the canonical dual problem is essentially equivalent to  the standard semi-definite programming (SDP) problem  and then be solved by well-developed software packages, such as SeDuMi \citep{sturn}. In the case that the canonical dual problem has no critical point,  four strategies are proposed to develop efficient algorithms based on unconstrained approaches by using the core points of canonical dual transformation methodology, complementary-dual principle and global optimality condition. A series of  fourth-order polynomials benchmark optimization problems are provided to demonstrate the effectiveness and efficiency of the proposed strategies.
\section{A brief review of canonical duality theory}
\label{sec:1}
For the completeness of this paper, we give a brief review of the following fourth-order polynomials minimization problem (primal problem) in \citet{gao3}:
\begin{eqnarray}
(\mathcal{P}): \min \Big\{  P(\mathbf{x}) = W(\mathbf{x}) + \frac{1}{2}\mathbf{x}^TQ\mathbf{x} -\mathbf{x}^T\mathbf{f}:\mathbf{x} \in \mathbb{R}^n\Big\},
\end{eqnarray}
where,
\begin{eqnarray}
W(\mathbf{x})= \sum_{k=1}^{m} \frac{1}{2} \alpha_k \Big(\frac{1}{2}\mathbf{x}^TA_k\mathbf{x} + \mathbf{b}^T_kx + c_k \Big)^2,
\end{eqnarray}
and $A_k = A^T_k, Q = Q^T \in \mathbb{R}^{n \times n}$ are indefinite symmetrical matrices, $\mathbf{b}_k,\mathbf{f} \in \mathbb{R}^n$ are given vectors, $\alpha_k, c_k \in \mathbb{R}$ are known constants. Without loss of much generality, the $\alpha_k$ is assumed to be positive.\\
\indent The standard canonical dual transformation methodology consists of the following three procedures.
\subsection{Canonical dual transformation}
\indent Introducing a nonlinear operator(a G\^{a}teaux differentiable geometrical measure)
\begin{eqnarray}
\bm{\xi} = (\xi_1,\cdots,\xi_m)^T = \mathrm{\Lambda}(\mathbf{x}) = \Big\{ \frac{1}{2} \mathbf{x}^T A_k \mathbf{x} + \mathbf{b}^T_k \mathbf{x} + c_k \Big\}^m: \mathbb{R}^n \rightarrow \mathcal{E}_a \subset \mathbb{R}^m
\end{eqnarray}
so that $W(\mathbf{x})$ can be recast by:
\begin{eqnarray}
W(\mathbf{x}) = V(\mathrm{\Lambda}(\mathbf{x})),
\end{eqnarray}
where, $V(\bm{\xi})$ is said to be a canonical function and in this case
\begin{eqnarray}
V(\bm{\xi}) = \sum_{k=1}^{m} \frac{1}{2} \alpha_k \xi^2_k = \frac{1}{2} \bm\alpha^T(\bm\xi \circ \bm\xi),
\end{eqnarray}
in which, $\bm\alpha = (\alpha_1, \cdots, \alpha_m)^T$, the notation $\mathbf{s} \circ \mathbf{t} = (s_1 s_1, \cdots, s_m s_m)^T$ denotes the Hadamard product for any two vectors
$\mathbf{s}$, $\mathbf{t} \in \mathbb{R}^m$.\\
\indent Then, the primal problem can be rewritten as the canonical form:
\begin{eqnarray}
\min_{\mathbf{x} \in \mathbb{R}^n} \Big\{  P(\mathbf{x}) = V(\mathrm{\Lambda}(\mathbf{x})) - U(\mathbf{x}) \Big\},
\end{eqnarray}
where $U(\mathbf{x}) = - \frac{1}{2}\mathbf{x}^TQ\mathbf{x} + \mathbf{x}^T\mathbf{f}$.
\subsection{Generalized complementary function}
\indent The dual variable $\bm\varsigma$ to $\bm\xi$ is defined by the duality mapping
\begin{eqnarray}
\bm \varsigma = (\varsigma_1, \cdots, \varsigma_m) = \nabla V(\bm{\xi}) = \bm\alpha \circ \bm\xi:\mathcal{E}_a \rightarrow \mathcal{E}^{*}_a \subset \mathbb{R}^m.
\end{eqnarray}
\indent For the given canonical function $V(\bm{\xi})$, the Legendre conjugate $V^{\ast}(\bm \varsigma)$ can be defined by:
\begin{eqnarray}
V^{\ast}(\bm{\varsigma}) = \underset{\bm{\xi}}{\mathrm{sta}}\{\bm{\xi}^T\bm{\varsigma} - V(\bm{\xi})\} = \sum_{k=1}^{m}\frac{1}{2}\alpha^{-1}_k \varsigma^2_k,
\end{eqnarray}
where, sta$\{\cdot\}$ stands for finding stationary point of the statement in $\{\cdot\}$. The $(\bm \xi,\bm \varsigma)$ forms a canonical duality pair and the following canonical
duality relations hold on $\mathcal{E}_a \times \mathcal{E}^{*}_a$:
\begin{eqnarray}
\bm{\varsigma}  = \nabla V(\bm{\xi}) \Leftrightarrow \bm{\xi}  = \nabla V^{\ast}(\bm{\varsigma}) \Leftrightarrow V(\bm{\xi}) + V^{\ast}(\bm{\varsigma}) = \bm{\xi}^T \bm{\varsigma}.
\end{eqnarray}
\indent Replacing $W(\mathbf{x})= V(\mathrm{\Lambda}(\mathbf{x}))$ by $\mathrm{\Lambda}^T(\mathbf{x})\bm{\varsigma} - V^{\ast}(\bm{\varsigma})$, the generalized complementary function can be defined by
\begin{eqnarray}
\Xi(\mathbf{x},\bm{\varsigma}) &=& \mathrm{\Lambda}^T(\mathbf{x})\bm{\varsigma} - V^{\ast}(\bm{\varsigma}) - U(\mathbf{x}) \nonumber \\
& = & \sum_{k=1}^{m} \Big[\Big(\frac{1}{2}\mathbf{x}^TA_k\mathbf{x} + \mathbf{b}^T_k \mathbf{x} + c_k \Big) \varsigma_k -\frac{1}{2}\alpha^{-1}_k \varsigma^2_k\Big] + \frac{1}{2}\mathbf{x}^TQ\mathbf{x} -\mathbf{x}^T\mathbf{f}. ~~
\end{eqnarray}
\subsection{Canonical dual function}
\indent By using the generalized complementary function, the canonical dual function $P^d(\bm{\varsigma})$ can be formulated as
\begin{eqnarray}
P^d(\bm{\varsigma}) = \underset{\mathbf{x}}{\mathrm{sta}}\{\Xi(\mathbf{x},\bm{\varsigma})\}.
\end{eqnarray}
\indent For a fixed $\bm{\varsigma}$, the stationary condition $\nabla\Xi(\mathbf{x},\bm{\varsigma})$ leads to the canonical equilibrium equation:
\begin{eqnarray}
G(\bm{\varsigma}) \mathbf{x} = F(\bm{\varsigma}),
\end{eqnarray}
in which, $G(\bm{\varsigma}) = Q + \sum_{k=1}^{m} \varsigma_k A_k$, $F(\bm{\varsigma}) = \mathbf{f} - \sum_{k=1}^{m} \varsigma_k \mathbf{b}_k$. For any given $\bm\varsigma$, if $F(\bm{\varsigma})$ is in the column space of $G(\bm{\varsigma})$, denoted by $\mathcal{C}_{ol}(G(\bm{\varsigma}))$, i.e., a linear space spanned by the columns of $G(\bm{\varsigma})$, the solution of the canonical equilibrium equation can be well defined by
\begin{eqnarray}
\mathbf{x} = G^{\dagger}(\bm{\varsigma})F(\bm{\varsigma}),
\end{eqnarray}
where, $G^{\dagger}(\bm{\varsigma})$ denotes the Moore-Penrose generalized inverse of $G(\bm{\varsigma})$.\\
\indent Then, the canonical dual function can be written explicitly as follows
\begin{eqnarray}
P^d(\bm{\varsigma}) = \sum_{k=1}^{m}\Big(c_k \varsigma_k - \frac{1}{2} \alpha^{-1}_k \varsigma^2_{k}\Big) - \frac{1}{2} F^T(\bm{\varsigma})G^{\dagger}(\bm{\varsigma})F(\bm{\varsigma}).
\end{eqnarray}
\indent Finally, the canonical dual problem can be expressed by
\begin{eqnarray}
(\mathcal{P}^d): \mathrm{sta} \Big\{P^d(\bm{\varsigma}) = \sum_{k=1}^{m}\Big(c_k \varsigma_k - \frac{1}{2} \alpha^{-1}_k \varsigma^2_{k}\Big) - \frac{1}{2} F^T(\bm{\varsigma})G^{\dagger}(\bm{\varsigma})F(\bm{\varsigma}): \bm{\varsigma} \in \mathcal{S}_a\Big\},~~~
\end{eqnarray}
where the dual feasible space is defined by $\mathcal{S}_a = \{ \bm{\varsigma} \in \mathbb{R}^m | F(\bm{\varsigma}) \in \mathcal{C}_{ol}(G(\bm{\varsigma}))\}$.\\

\noindent \textbf{Theorem 1} (\textbf{Complementary-Dual Principle and Analytical Solution}).
The problem $(\mathcal{P}^d)$ is canonically dual to the primal problem $(\mathcal{P})$ in the sense that if $\bar{\bm{\varsigma}}$ is a critical
point of $(\mathcal{P}^d)$, then the vector
\begin{equation}
\bar{\mathbf{x}} = G^{\dagger}(\bar{\bm{\varsigma}})F(\bar{\bm{\varsigma}})
 \end{equation}
 is a critical point of $(\mathcal{P})$ and
\begin{equation}
 P(\bar{\mathbf{x}})=P^d(\bar{\bm{\varsigma}}) .
 \end{equation}

 This theorem shows that the critical solutions to the primal problem depend  analytically on the
 canonical dual solutions and there is no duality gap between the primal problem and its canonical dual.\\

\noindent \textbf{Theorem 2} (\textbf{Global Optimality Condition}). Suppose $  \bar{\bm{\varsigma}} $ is a critical point of
$P^d( {\bm{\varsigma}})$. If $\bar{\bm{\varsigma}} \in S^{+}_a$,
then $\bar{\bm{\varsigma}}$ is a global maximizer of $(\mathcal{P}^d)$ on $S^{+}_a$ if and only if
the analytical solution $\bar{\mathbf{x}} = G^{\dagger}(\bar{\bm{\varsigma}})F(\bar{\bm{\varsigma}}) $ is a global minimizer of $(\mathcal{P})$ on $\mathbb{R}^n$, i.e.,
\begin{equation}
P(\mathbf{\bar{x}}) = \min_{\mathbf{x} \in \mathbb{R}^n}P(\mathbf{x}) \Leftrightarrow \max_{\mathbf{\bm{\varsigma}} \in \mathcal{S}^+_a}P^d(\bm{\varsigma}) = P^d(\mathbf{\bar{\bm{\varsigma}}}),
 \end{equation}
where
\begin{equation}
\mathcal{S}^{+}_a = \{\bm \varsigma \in \mathcal{S}_a | G(\bm \varsigma) \succeq 0 \}.
 \end{equation}

This theorem shows that $\bar{\bm{\varsigma}} \in S^{+}_a$ provides a global optimality condition, which can be used to develop algorithms for solving the nonconvex primal problem.

\section{The equivalent semi-definite programming problem}
By Theorem 2, the primal problem is equivalent to the following
  canonical dual maximization problem ($({\cal P}^d_{\max})$ in short):
\begin{eqnarray}
({\cal P}^d_{\max}): \max_{\mathbf{\bm{\varsigma}} \in \mathcal{S}^+_a}P^d(\bm{\varsigma}) = \sum_{k=1}^{m}\Big(c_k \varsigma_k - \frac{1}{2} \alpha^{-1}_k \varsigma^2_{k}\Big) - \frac{1}{2} F^T(\bm{\varsigma})G^{\dagger}(\bm{\varsigma})F(\bm{\varsigma})
\end{eqnarray}
In this section,  we will show that this problem
can be also equivalent to the standard semi-definite programming problem (SDP).\\
\begin{eqnarray}
(SDP): \min && \frac{1}{2}t_1 + \frac{1}{2}t_2 - \bm\varsigma^T \mathbf{c}  \nonumber \\
\mathrm{subject~to:} &&
\begin{pmatrix}
    G(\bm{\varsigma}) & F(\bm{\varsigma}) \nonumber  \\
    F^T(\bm{\varsigma}) & t_1
\end{pmatrix} \succeq 0, \\
&&
\begin{pmatrix}
    Diag\{\alpha_1,\cdots,\alpha_m\} & \bm{\varsigma} \\
    \bm{\varsigma}^T & t_2
\end{pmatrix} \succeq 0.
\end{eqnarray}
\noindent \textbf{Theorem 3} Let $(\bar{\bm{\varsigma}}, \bar{t}_1, \bar{t}_2)$ be an optimal solution of problem (SDP), if $G(\bar{\bm{\varsigma}}) \succ 0$, then $\bar{\bm{\varsigma}}$ is the unique optimal solution of problem $(\mathcal{P}^d)$ and $\bar{\mathbf{x}} = G^{\dagger}(\bar{\bm{\varsigma}})F(\bar{\bm{\varsigma}})$ is the unique optimal solution of problem $(\mathcal{P})$. If det($G(\bar{\bm{\varsigma}})) = 0$ and $(I - G(\bar{\bm{\varsigma}})G^{\dagger}(\bar{\bm{\varsigma}}))F(\bar{\bm{\varsigma}}) = 0$, then $\bar{\bm{\varsigma}}$ is an optimal solution of problem $(\mathcal{P}^d)$ and $\bar{\mathbf{x}} = G^{\dagger}(\bar{\bm{\varsigma}})F(\bar{\bm{\varsigma}})$ is an optimal solution of problem $(\mathcal{P})$. In this case, problem $(\mathcal{P})$ has multiple optimal solutions.\\

\noindent \textbf{Proof}. At first, we relax $({\cal P}^d_{\max})$ to the following form
\begin{eqnarray}
\min && \frac{1}{2}t_1 + \frac{1}{2}t_2 - \bm\varsigma^T \mathbf{c}  \nonumber \\
\mathrm{subject~to:} && t_1 \geq F^T(\bm{\varsigma})G^{-1}(\bm{\varsigma})F(\bm{\varsigma}),  \nonumber  \\
                     && t_2 \geq  \bm{\varsigma}^T Diag\{\alpha_1,\cdots,\alpha_m\} \bm{\varsigma}, \nonumber  \\
                     && G(\bm{\sigma}) \succeq 0 ,
\end{eqnarray}
where, $\mathbf{c} = [c_1, \dots, c_m]^T$, $Diag\{\alpha_1,\cdots,\alpha_m\}$ stands for a diagonal matrix with $\alpha_1,\cdots,\alpha_m$ as its elements.\\
\indent \textbf{Lemma 1} (Schur complement) Considering the partitioned symmetric matrix
\begin{eqnarray}
X = X^T =
\begin{pmatrix}
A & B\\
B^T & C \\
\end{pmatrix},
\end{eqnarray}
if $A \succ 0$, then $X\succeq 0$  if and only if the matrix $C-B^TA^{-1}B\succeq 0 $.
\\
\indent Using the Schur complement lemma, we can get the equivalent positive (semi) definite programming optimization problem (SDP) consequently according to Theorem 1 and Theorem 2.
~~~~~~~~~~~~~~~~~~~~~~~~~~~~~~~~~~~~~~~~~~~~~~~~~~~~~~~~~~~~~~~~~~~~~~~~~~~~~~~~~~~~~~~~~~~~~~~~~$\Box$

\section{Four strategies for canonical dual theory}
Although the canonical dual problem can be transformed into the equivalent semi-definite programming problem and then solved by well-developed software packages, it should be noted that there may be no critical points in the canonical dual feasible domain. Moreover, solving a constrained optimization problem (SDP) is more complicated than an unconstrained one. In this paper, we focus on unconstrained methods, trying to explore efficient and effective algorithms based on the canonical duality theory.\\
\indent The canonical duality theory has provided with a unified analytic solution for optimization problems. Obviously, we can firstly find all of the stationary points, and then identify which one is in the canonical dual feasible domain $\mathcal{S}^{+}_a$. As can be seen from the main procedures of canonical dual transformation methodology, we can find that we have to solve stationary problems twice, one is for $\Xi(\mathbf{x},\bm{\varsigma})$, and the other is for $P^d(\bm{\varsigma})$. To calculate the stationary points for $\Xi(\mathbf{x},\bm \varsigma)$, we have to solve the following nonlinear equations:
\begin{equation}{\label{eq:16}}
\left\{
\begin{aligned}
&G(\bm{\varsigma}) \mathbf{x} = F(\bm{\varsigma}),\\
&\frac{1}{2}\mathbf{x}^TA_k\mathbf{x} + \mathbf{b}^T_k \mathbf{x} + c_k = \alpha^{-1}_k \varsigma_k,  \forall k = 1,\cdots,m.
\end{aligned}
\right.
\end{equation}
Similarly, to find the stationary points for $P^d(\bm{\varsigma})$, we have to solve the nonlinear equations as follows:
\begin{equation}{\label{eq:17}}
\frac{1}{2} F^T(\bm{\varsigma})G^{\dagger} A_k G^{\dagger} (\bm{\varsigma}) F(\bm{\varsigma}) + \mathbf{b}^T_k G^{\dagger} (\bm{\varsigma})F(\bm{\varsigma})  + c_k - \alpha^{-1}_k \varsigma_k =0, \forall k = 1,\cdots,m.
\end{equation}
\indent  As the number of variables becomes large, the complexity of computing the stationary problems is also increasing. That is to say, the computing of the stationary points for $\Xi(\mathbf{x},\bm{\varsigma})$ will be more complicated than that of $P^d(\bm{\varsigma})$. On the contrary, we can also observe that
the solving of stationary points for $P^d(\bm{\varsigma})$ may become more difficult than $\Xi(\mathbf{x},\bm{\varsigma})$ due to the inverse of matrix. It indicates that the complexity of solving the two nonlinear equations will be distinctive for different problems, which is the original source of why we design \textit{Strategy 1} and \textit{Strategy 2}.\\
\indent Any way, finding stationary points is just one way to solve optimization problems, and we can use iterative method to ``search" for global optimum as well. According to the results of canonical duality theory, compared with the primal problem $(\mathcal{P})$, the advantages of solving the canonical dual problem $(\mathcal{P}^d)$ is that $(\mathcal{P}^d)$ can be easily solved by well-developed optimization algorithms because the $(-\mathcal{P}^d)$ is convex on convex domain $\mathcal{S}^{+}_a$. In practice, we find that sometimes the form of $(\mathcal{P}^d)$ may become much more complicated than the primal problem $(\mathcal{P})$, also due to complexity of computing $G^{\dagger}(\bm{\varsigma})$, which is the original source of why we design \textit{Strategy 3} and \textit{Strategy 4}.\\
\indent So far, there still exist big issues in practical application. In terms of \textit{Strategy 1} and \textit{Strategy 2}, there may exist numerous stationary points, leading it difficult to identify which one is in the canonical dual feasible domain, while for \textit{Strategy 3}, there will be no canonical dual feasible solutions, making it impossible to substitute back to get solution to the primal problem. If we suppose that there exists a algorithm, once it runs into the ``neighborhood" of the canonical dual feasible domain, it will never deviate too far from the ``neighborhood", then we can start from an initial point in this ``neighborhood" to finally arrive at the global solution according to the close relationship between the solution of dual and that of the primal by the complementary-dual principle, which is the kernel of the proposed four strategies. In this case, there is no need to identify which stationary point is in the canonical dual feasible domain $\mathcal{S}^{+}_a$ and there is also no need to solving a constrained optimization problem, because we can just start from a ``good" initial point to solve the nonlinear equations or to find global optimum based on unconstrained approach.\\
\indent The detailed strategies of how to use the canonical duality theory above to find a global minimum will be given in the following:
\\
\begin{table}[!htbp]
\begin{tabular}{l}
\hline
\textbf{Strategy 1}\\
\hline
1: \textit{standardization}\\
\indent Convert the original problem to the standard form of the primal problem discussed in \\the paper,
and then the parameters of $(\mathcal{P})$ like $\alpha_k, A_k, \mathbf{b}_k, c_k, Q, \mathbf{f}$ will be well defined.\\
2: \textit{selection}\\
\indent Select an appropriate $\bm \varsigma_0$ to make sure that $G(\bm \varsigma_0) \succeq 0$, and then gain the corresponding\\
$\mathbf{x}_0 = G^{-1}(\bm \varsigma_0) F(\bm \varsigma_0)$.\\
3: \textit{nonlinear equations}\\
\indent Taking $(\mathbf{x}_0,\bm \varsigma_0)$ as initial point, use numerical calculation methods to solve
the \\nonlinear equations in (\ref{eq:16}).\\
\hline
\end{tabular}
\end{table}
\\
\begin{table}[!htbp]
\begin{tabular}{l}
\hline
\textbf{Strategy 2}\\
1: \textit{standardization}\\
\indent Convert the original problem to the standard form of the primal problem discussed in \\the paper, and then the parameters of $(\mathcal{P})$ like $\alpha_k, A_k, \mathbf{b}_k, c_k, Q, \mathbf{f}$ will be well defined.\\
2: \textit{selection}\\
\indent Select an appropriate $\bm \varsigma_0$ to make sure that $G(\bm \varsigma_0) \succeq 0$.\\
3: \textit{nonlinear equations}\\
\indent Taking $\bm\varsigma_0$ as initial point, use numerical calculation methods to solve the \\nonlinear equations in (\ref{eq:17}).\\
\hline
\end{tabular}
\end{table}
\\
\begin{table}[!htbp]
\begin{tabular}{l}
\hline
\textbf{Strategy 3}\\
1: \textit{standardization}\\
\indent Convert the original problem to the standard form of the primal problem discussed in \\the paper, and then the parameters of $(\mathcal{P})$ like $\alpha_k, A_k, \mathbf{b}_k, c_k, Q, \mathbf{f}$ will be well defined.\\
2: \textit{selection}\\
\indent Select an appropriate $\bm \varsigma_0$ to make sure that $G(\bm \varsigma_0) \succeq 0$.\\
3: \textit{numerical optimization}\\
\indent Taking $\bm \varsigma_0$ as an initial point (or some $\bm \varsigma_0$ as initial population), using numerical\\ optimization algorithms to optimize the dual problem. \\
\hline
\end{tabular}
\end{table}
\\
\begin{table}[!htbp]
\begin{tabular}{l}
\hline
\indent \textbf{Strategy 4}\\
1: \textit{standardization}\\
\indent Convert the original problem to the standard form of the primal problem discussed in \\the paper, and then the parameters of $(\mathcal{P})$ like $\alpha_k, A_k, \mathbf{b}_k, c_k, Q, \mathbf{f}$ will be well defined.\\
2: \textit{selection}\\
\indent Select an appropriate $\bm \varsigma_0$ to make sure that $G(\bm \varsigma_0) \succeq 0$, and then gain the corresponding \\
$\mathbf{x}_0 = G^{-1}(\bm \varsigma_0) F(\bm \varsigma_0)$.\\
3: \textit{numerical optimization}\\
\indent Taking $\mathbf{x}_0$ as an initial point (or some $\mathbf{x}_0$ as initial population), using numerical \\optimization algorithms to optimize the primal problem. \\
\hline
\end{tabular}
\end{table}
\\
\indent \textbf{Remark 1.} In the numerical optimization step of \textit{Strategy 3} and \textit{Strategy 4}, some numerical optimization methods may not be able to search just in the ``neighborhood", in this case, a penalty function is suggested to add to $(\mathcal{-P}^d)$ or $(\mathcal{P})$ so that $G(\bm \varsigma) \succeq 0$ in the search process.
\section{Numerical results}
\begin{figure}[!htbp]
  \includegraphics[width=0.4\textwidth,height=0.3\textwidth]{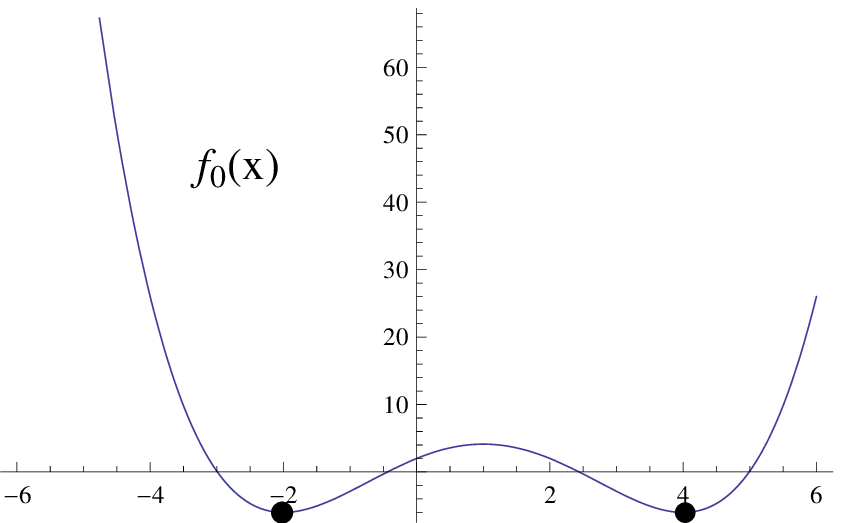}
  \includegraphics[width=0.4\textwidth,height=0.3\textwidth]{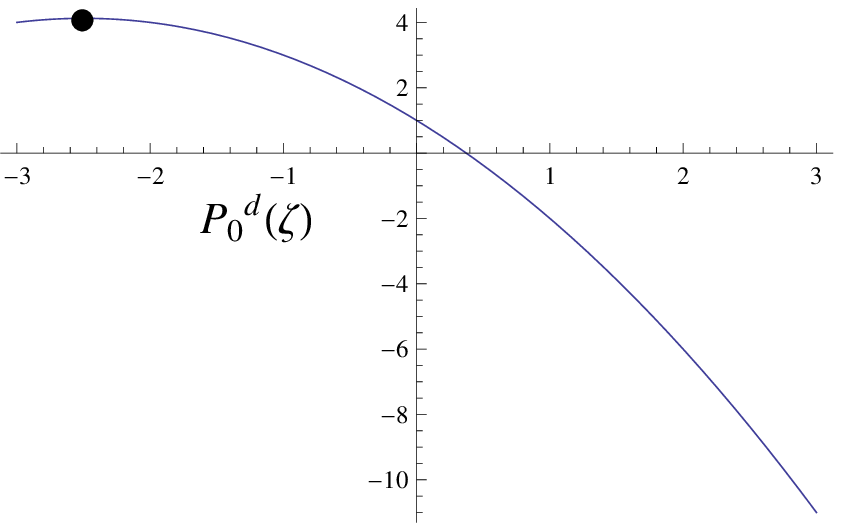}
\caption{graphs of the primal function $f_0(x)$ and its corresponding dual function $P^{d}_0(\varsigma)$}
\label{fig:f0}       % Give a unique label
\end{figure}

\begin{figure}[!htbp]
  \includegraphics[width=0.4\textwidth,height=0.3\textwidth]{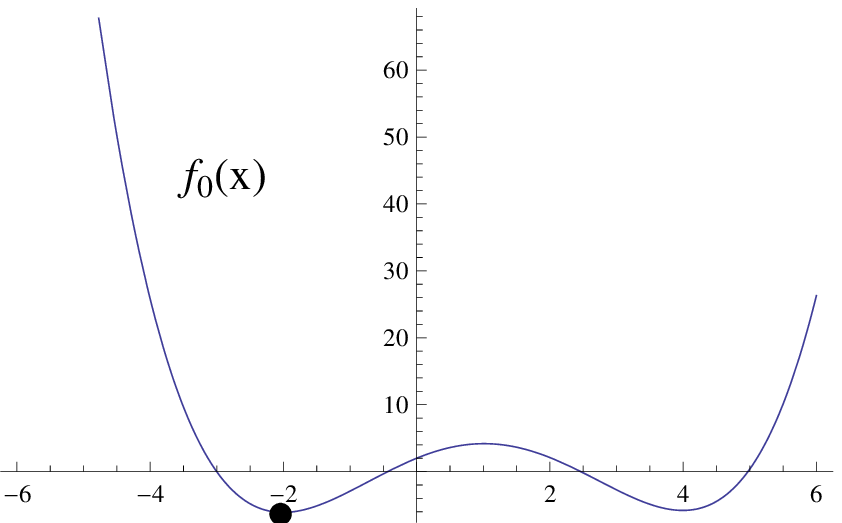}
  \includegraphics[width=0.4\textwidth,height=0.3\textwidth]{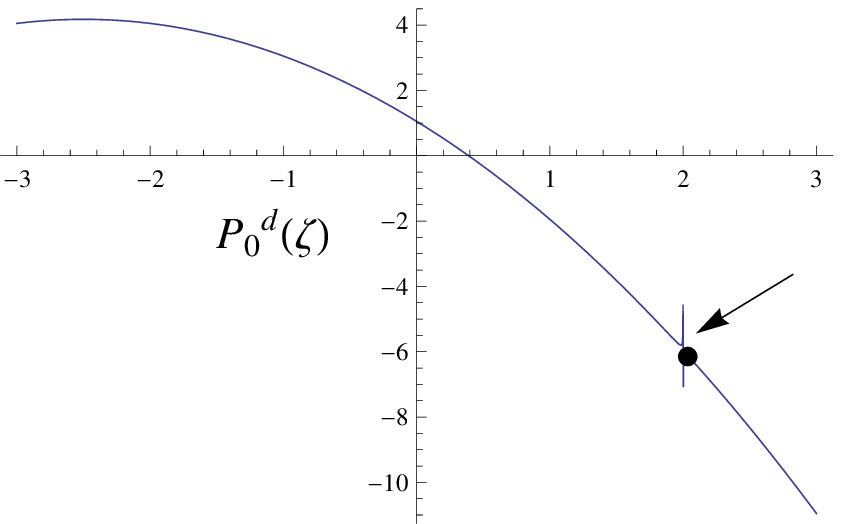}
\caption{graphs of the primal function $f_0(x)$ and its dual function $P^{d}_0(\varsigma)$ with linear perturbation}
\label{fig:f0p}       % Give a unique label
\end{figure}

\begin{figure}[!htbp]
  \includegraphics[width=0.8\textwidth,height=0.6\textwidth]{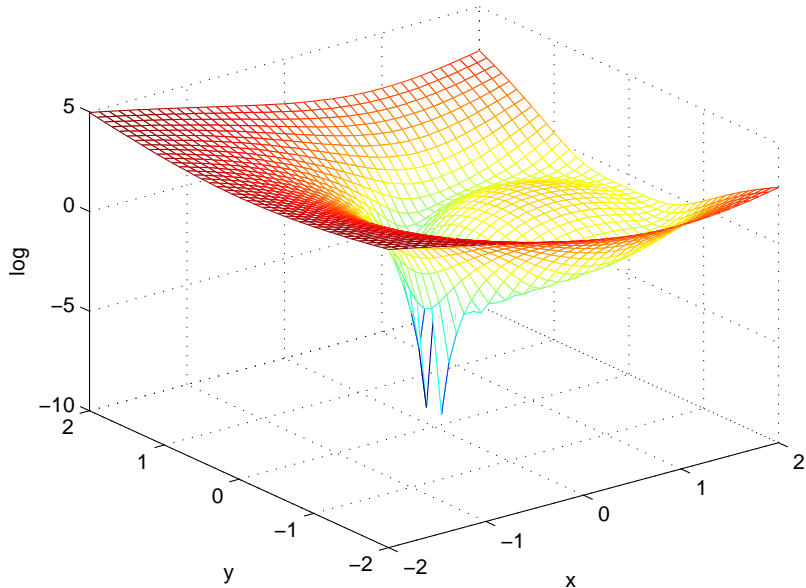}
\caption{Landscape of Zettle function}
\label{fig:zettle}       % Give a unique label
\end{figure}

\begin{figure}[!htbp]
  \includegraphics[width=0.8\textwidth,height=0.6\textwidth]{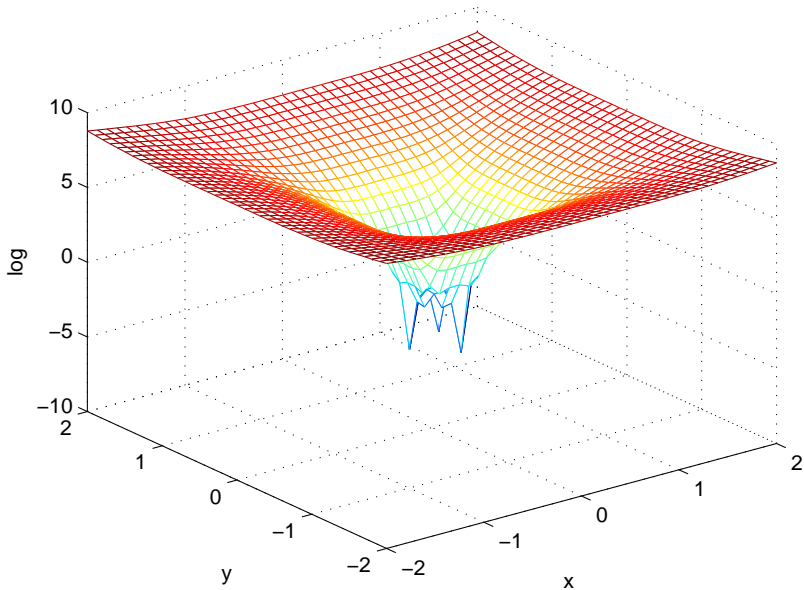}
\caption{Landscape of Styblinski-Tang function}
\label{fig:st}       % Give a unique label
\end{figure}

\begin{figure}[!htbp]
  \includegraphics[width=0.8\textwidth,height=0.6\textwidth]{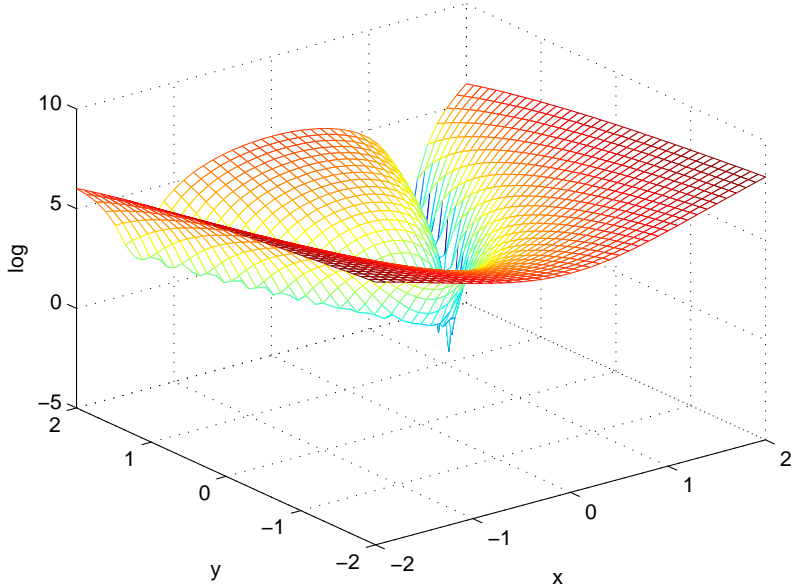}
\caption{Landscape of Rosenbrock function in two dimension}
\label{fig:1}       % Give a unique label
\end{figure}

\begin{figure}[!htbp]
  \includegraphics[width=0.8\textwidth,height=0.6\textwidth]{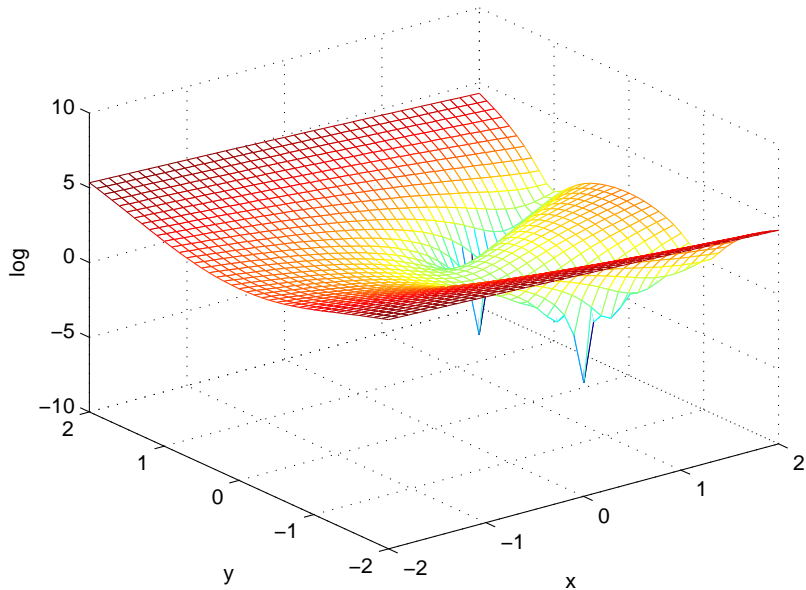}
\caption{Landscape of Dixon and Price function in two dimension}
\label{fig:dp}       % Give a unique label
\end{figure}

To testify the effectiveness of the strategies, some fourth-order polynomials benchmark functions are collected, and we will use the proposed strategies to find the global
minimum one by one. In this paper, we implement the strategies in MATLAB R2010b on Intel(R) Core(TM) i3-2310M CPU @2.10GHz under Window 7 environment, and \textit{fsolve} and \textit{fminunc} built in MATLAB are used to solve nonlinear equations and for numerical optimization, respectively.\\
\indent \textit{Example 0} (A special case)\\
Considering the following one-dimensional problem
\begin{eqnarray}
f_0(x)&& = \frac{1}{2}\alpha(\frac{1}{2}a x^2 + b x + c)^2 + \frac{1}{2}q x^2 - x f. \nonumber
\end{eqnarray}
we can get the corresponding canonical dual problem easily
\begin{eqnarray}
P^{d}_0(\varsigma)&& = c\varsigma - \frac{1}{2\alpha}\varsigma^2 - \frac{(f - b\varsigma)^2}{2(q+a\varsigma)}.\nonumber
\end{eqnarray}
To be more specific, let fix $\alpha = 1$, $a = 1$, $b = -1$, $c = -2$, $q = -2$, $f = -2$, which is a special case because $x^{*} = G^{-1}(\varsigma)F(\varsigma) = \frac{f-b\varsigma}{q+a\varsigma}=1, \forall \varsigma$,
and the graphs of the primal and dual functions are given in Fig.\ref{fig:f0}.
\\
\indent As can be shown in Fig.\ref{fig:f0}, there is no critical point in the canonical dual feasible domain $G(\varsigma) = -2 + \varsigma > 0$, that is to say, the semi-definite programming software packages will be invalid in this case. For remedy, we can add a small linear perturbation to the primal problem, for instance, $f = f + \triangle $. If $\triangle = 0.05$, the graphs of the primal and dual functions with linear perturbation are given in Fig.\ref{fig:f0p}.
\\
\indent Using SeDuMi to solve the modified canonical dual problem, we can get $\varsigma^{*} = 2.0166$, $x^{*} = G^{-1}(\varsigma)F(\varsigma) = -2.0056$ and $P(x^{*}) = -6.1234$. We can find that there exists small deviation from the global optimum. \\
\indent On the other hand, if we choose to use the proposed \textit{Strategy 4}, not starting from the initial point $x_0 = 1$ directly but adding a translation $x_0 = 1 + randn$ (\textit{randn} is the built-in function of the standard normal distribution within MATLAB), we can finally arrive at either $x^{*} = -2$ or $x^{*} = -4$ precisely.
\\
\\
\indent \textit{Example 1} (Colville function)
\begin{eqnarray}
f_1&& = 100(x_2 - x^2_1)^2 + (1 - x_1)^2 + 90(x_4 - x^2_3)^2 + (1 - x_3)^2 \nonumber \\
&&+ 10.1((x_2 - 1)^2 + (x_4 - 1)^2) + 19.8(x_2 -1)(x_4 -1). \nonumber
\end{eqnarray}
\indent We firstly rewrite it to the standard form, and then we can get $\alpha_1 = 200, \alpha_2 = 180, A_1 = Diag\{-2,0,0,0\}, A_2 = Diag\{0,0,-2,0\}, \mathbf{b}_1 = [0,1,0,0]^T$,
$\mathbf{b}_2 = [0,0,0,1]^T, c_1 = c_2 =0,\mathbf{f} = [2,40,2,40]^T$, and
%$Q = [2,0,0,0;0,20.2,0,19.8;0,0,2,0; 0;19.8;0;20.2]$,
\begin{eqnarray}
Q =
\begin{pmatrix}
2&0&0&0\\
0&20.2&0&19.8\\
0&0&2&0\\
0&19.8&0&20.2
\end{pmatrix}, \nonumber
\mathrm{then} \quad
G(\bm \varsigma) =
\begin{pmatrix}
2 -2\varsigma_1 & 0  & 0 & 0   \\
0 &20.2& 0 &19.8 \\
0 &0   & 2 - 2\varsigma_2 &0    \\
0 &19.8& 0 &20.2
\end{pmatrix}. \nonumber
\end{eqnarray}
\indent \textit{Strategy 1}\\
\indent The generalized complementary function is
\begin{eqnarray}
&& ~~\Xi(\mathbf{x},\bm \varsigma) = (x^2_1 - x_2)\varsigma_1  + (x^2_3 - x_4)\varsigma_2 - \frac{1}{400}\varsigma^2_1 - \frac{1}{360}\varsigma^2_2 \nonumber \\
&&+ (x^2_1 + 10.1x^2_2 + x^2_3 + 10.1x^2_4 + 19.8x_2x_4) - (2x_1 + 40x_2 + 2x_3 + 40x_4) + 42.
\nonumber
\end{eqnarray}
\indent We select $\bm \varsigma_0 = (\varsigma_1,\varsigma_2) = (0.5,0.5)$ to make sure that $G(\bm \varsigma_0) \succeq 0$ and the corresponding $\mathbf{x}_0 = (x_1,x_2,x_3,x_4) = G(\bm \varsigma_0)^{-1} F(\bm \varsigma_0)= (2.0000,0.9875,2.0000,0.9875)$  as initial point $(\mathbf{x}_0, \bm \varsigma_0)$ for the nonlinear equations in (\ref{eq:16}), and after 5 iterations with 0.338344 seconds, we obtain $(\mathbf{x}^{*}, \bm \varsigma^{*}) = (1.0000,1.0000,1.0000,1.0000,0.0000,0.0000)$ and $P(\mathbf{x}^\ast) = 0$.
\\
\indent \textit{Strategy 2}\\
\indent The canonical dual function is
\begin{eqnarray}
P^d(\bm \varsigma) = &&  42 -  \frac{1}{400}\varsigma^2_1 - \frac{1}{360}\varsigma^2_2 \nonumber \\
&& - \frac{1}{2} \Big(2, 40 - \varsigma_1, 2, 40 - \varsigma_2\Big)
\begin{pmatrix}
2 -2\varsigma_1 & 0  & 0 & 0   \\
0 &20.2& 0 &19.8 \\
0 &0   & 2 - 2\varsigma_2 &0    \\
0 &19.8& 0 &20.2
\end{pmatrix}^{+}
\begin{pmatrix}
2\\
40 - \varsigma_1\\
2\\
40 - \varsigma_2
\end{pmatrix}.
\nonumber
\end{eqnarray}
\indent We select the same $\bm \varsigma_0 = (0.5,0.5)$ for the nonlinear equations in (\ref{eq:17}), and after 6 iterations with 0.302655 seconds, we obtain $\bm \varsigma^{*} = (0.0000,0.0000)$. The corresponding $\mathbf{x}^\ast = G(\bm \varsigma^\ast)^{-1} F(\bm \varsigma^\ast)= (1,1,1,1)$ and $P(\mathbf{x}^\ast) = 0$.
\\
\indent \textit{Strategy 3}\\
\indent We choose the same $\bm \varsigma_0 = (0.5,0.5)$ as initial point for $(-\mathcal{P}^d)$, and we can finally arrive at
$\bm \varsigma^\ast = (0,0)$ with 9 iterations and 0.329829 seconds. The corresponding $\mathbf{x}^\ast = G(\bm \varsigma^\ast)^{-1} F(\bm \varsigma^\ast)= (1,1,1,1)$ and $P(\mathbf{x}^\ast) = 0$.\\
\indent \textit{Strategy 4}\\
\indent We choose the same $\bm \varsigma_0 = (0.5,0.5)$, and then we use the corresponding $\mathbf{x}_0 = G(\bm \varsigma_0)^{-1} F(\bm \varsigma_0)= (2.0000,0.9875,2.0000,0.9875)$ as initial point for $P(\mathbf{x})$, and we can finally arrive at $\mathbf{x}^\ast = (1,1,1,1)$ and $P(\mathbf{x}^\ast) = 0$ with 26 iterations and 0.354057 seconds.\\
\\
\indent \textit{Example 2} (Zettle function)
\begin{eqnarray}
f_2  = (x^2_1 + x^2_2 - 2x_1)^2 + 0.25x_1. \nonumber
\end{eqnarray}
\indent The landscape of Zettle function is given in Fig. \ref{fig:zettle}.
\\
\indent Firstly, we rewrite it to the standard form, and then we can get
$\alpha = 2, A = Diag\{2,2\}, b = [-2,0]^T, c = 0, Q = \mathbf{0}, f = [-0.25,0]^T$, and then $G(\bm \varsigma) = Diag\{2\varsigma_1,2\varsigma_1\}$.\\
\indent \textit{Strategy 1}\\
\indent The generalized complementary function is
\begin{eqnarray}
\Xi(\mathbf{x},\bm \varsigma) = (\frac{1}{2} (2 x_1^2+2 x_2^2)-2 x_1) \varsigma -\frac{\varsigma^2}{4} + 0.25 x_1. \nonumber
\end{eqnarray}
\indent We select $\bm \varsigma_0 = 0.1$ to make sure that $G(\bm \varsigma_0) \succeq 0$ and the corresponding $\mathbf{x}_0 = G(\bm \varsigma_0)^{-1} F(\bm \varsigma_0)$ $=(-0.2500,0)$ as initial point $(\mathbf{x}_0, \bm \varsigma_0)$ for the nonlinear equations in (\ref{eq:16}), and after 3 iterations with 0.290396 seconds, we obtain $(\mathbf{x}^{*}, \bm \varsigma^{*}) = (-0.0299,0,0.1214)$ and $P(\mathbf{x}^\ast) = -0.0038$.
\\
\indent \textit{Strategy 2}\\
\indent The canonical dual function is
\begin{eqnarray}
P^d(\bm \varsigma) = -\frac{\varsigma^2}{4}-\frac{(2 \varsigma-0.25)^2}{4 \varsigma}. \nonumber
\end{eqnarray}
\indent We select the same $\bm \varsigma_0 = 0.1$ for the nonlinear equations in (\ref{eq:17}), and after 4 iterations with 0.290140 seconds, we obtain $\bm \varsigma^{*} = 0.1214$. The corresponding $\mathbf{x}^\ast = G(\bm \varsigma^\ast)^{-1} F(\bm \varsigma^\ast)= -0.0299$ and $P(\mathbf{x}^\ast) = -0.0038$.
\\
\indent \textit{Strategy 3}\\
\indent We choose the same $\bm \varsigma_0 = 0.1$ as initial point for $(-\mathcal{P}^d)$, and we can finally arrive at
$\bm \varsigma^\ast = 0.1214$  with 4 iterations and 0.319423 seconds. The corresponding
$\mathbf{x}^\ast = G(\bm \varsigma^\ast)^{-1} F(\bm \varsigma^\ast)= (-0.0299,0)$ and $P(\mathbf{x}^\ast) = -0.0038$.\\
\indent \textit{Strategy 4}\\
\indent We choose the same $\bm \varsigma_0 = 0.1$, and then we use the corresponding $\mathbf{x}_0 = G(\bm \varsigma_0)^{-1} F(\bm \varsigma_0)$ $=(-0.2500,0)$ as initial point for $P(\mathbf{x})$, and we can finally arrive at $\mathbf{x}^\ast = (-0.0299,0)$ and $P(\mathbf{x}^\ast) = -0.0038$ with 5 iterations and 0.309056 seconds.\\
\\
\indent \textit{Example 3} (Styblinski-Tang function)
\begin{eqnarray}
f_3 = \frac{1}{2}\sum_{i=1}^{2}(x^4_i - 16x^2_i + 5x_i). \nonumber
\end{eqnarray}
\indent The landscape of Styblinski-Tang function is given in Fig. \ref{fig:st}.
\\
\indent At first, we rewrite it to the standard form, and then we can get
$\alpha_1 = \alpha_2 = 1, A_1 = Diag\{2,0\}, A_2 = Diag\{0,2\}, b_1 = b_2 = \mathbf{0} , c_1 = c_2 = 0, Q = Diag\{-16,-16\}, f = [-2.5,-2.5]^T$,
and then $G(\bm \varsigma) = Diag\{-16+2\varsigma_1,-16+2\varsigma_2\}$.\\
\indent \textit{Strategy 1}\\
\indent The generalized complementary function is
\begin{eqnarray}
\Xi(\mathbf{x},\bm \varsigma) = x_1^2 \varsigma_1+x_2^2 \varsigma_2-\frac{\varsigma_1^2}{2}-\frac{\varsigma_2^2}{2}+\frac{1}{2} \left(-16 x_1^2-16 x_2^2\right) + 2.5 x_1+2.5 x_2.\nonumber
\end{eqnarray}
\indent We select $\bm \varsigma_0 = (8.1,8.1)$ to make sure that $G(\bm \varsigma_0) \succeq 0$ and the corresponding $\mathbf{x}_0 = G(\bm \varsigma_0)^{-1} F(\bm \varsigma_0)$ $=(-12.5,-12.5)$ as initial point $(\mathbf{x}_0, \bm \varsigma_0)$ for the nonlinear equations in (\ref{eq:16}), and after 8 iterations with 0.305630 seconds, we obtain $(\mathbf{x}^{*}, \bm \varsigma^{*}) = (-2.9035,-2.9035,8.4305,8.4305)$ and $P(\mathbf{x}^\ast) = -78.3323$.
\\
\indent \textit{Strategy 2}\\
\indent The canonical dual function is
\begin{eqnarray}
P^d(\bm \varsigma) = -\frac{\varsigma_1^2}{2}-\frac{\varsigma_2^2}{2} - \frac{1}{2} \left(\frac{6.25 \left(-16+2 \varsigma_1\right)}{256-32 \varsigma_1-32 \varsigma_2+4 \varsigma_1 \varsigma_2}+\frac{6.25 \left(-16+2 \varsigma_2\right)}{256-32 \varsigma_1-32 \varsigma_2+4 \varsigma_1 \varsigma_2}\right) \nonumber
\end{eqnarray}
\indent We select the same $\bm \varsigma_0 = (8.1,8.1)$ for the nonlinear equations in (\ref{eq:17}), and after 8 iterations with 0.305489 seconds, we obtain $\bm \varsigma^{*} = (8.4305,8.4305)$. The corresponding $\mathbf{x}^\ast = G(\bm \varsigma^\ast)^{-1} F(\bm \varsigma^\ast)= (-2.9035,-2.9035)$ and $P(\mathbf{x}^\ast) = -78.3323$.
\\
\indent \textit{Strategy 3}\\
\indent We choose the same $\bm \varsigma_0 = (8.1,8.1)$ as initial point for $(-\mathcal{P}^d)$, and then we can finally arrive at $(\bm\varsigma^{\ast}) = (8.4305,8.4305)$ within
7 iterations and 0.320454 seconds. The corresponding $\mathbf{x}^\ast = G(\bm \varsigma^\ast)^{-1} F(\bm \varsigma^\ast)= (-2.9035,-2.9035)$ and $P(\mathbf{x}^\ast) = -78.3323$.\\
\indent \textit{Strategy 4}\\
\indent We choose the same $\bm \varsigma_0 = (8.1,8.1)$, and then we use the corresponding $\mathbf{x}_0 = G(\bm \varsigma_0)^{-1} F(\bm \varsigma_0) = (-12.5,-12.5)$ as initial point for $P(\mathbf{x})$, we can finally arrive at $\mathbf{x}^\ast = (-2.9035, -2.9035)$ and $P(\mathbf{x}^\ast) = -78.3323$ with 10 iterations and 0.323847 seconds.\\
\\
\indent \textit{Example 4} (Rosenbrock function)
\begin{eqnarray}
f_4 = \sum_{i=1}^{n-1}[100(x_{i+1} - x_i^2)^2 + (x_i - 1)^2]. \nonumber
\end{eqnarray}
\indent At first, we rewrite it to the standard form, and then we can get $\alpha_k = 200, A_k = -2I_{k}, b_k = e_{k+1}, c_k = 0, Q = Diag\{\underbrace{2,2,\cdots,2}_{n-1},0\}, f = [\underbrace{2,2,\cdots,2}_{n-1},0]^T$, where $k = 1,2,\cdots,n-1$, $I_{k} \in \mathbb{R}^{n \times n}$ is a diagonal matrix with all zeros except the position $(k,k)$ having value 1 and $e_{k} \in \mathbb{R}^{n}$ is a unit vector with all zeros except the position $k$ having value 1. Then we can obtain
$G(\bm \varsigma) = Q + \sum_{k=1}^{n-1} \varsigma_k A_k = Diag\{2 - 2\varsigma_1,2 - 2\varsigma_2,\cdots, 2-2\varsigma_{n-1},0\}$,
$F(\bm \varsigma) = f - \sum_{k=1}^{n-1} \varsigma_k b_k = [2,2-\varsigma_1,\cdots,2-\varsigma_{n-2},-\varsigma_{n-1}]$.\\
\indent Without much loss of generality, $ n = 2$ is chosen for simple study, and its corresponding landscape is plotted
in Fig.\ref{fig:1}, in which, the global minimum is located in a long, deep, narrow, banana shaped flat valley.
\\
\indent \textit{Strategy 1}\\
\indent The generalized complementary function is
\begin{eqnarray}
\Xi(\mathbf{x},\varsigma) = \left(x_2-x_1^2\right) \varsigma-\frac{\varsigma^2}{400}+x_1^2-2 x_1 + 1 \nonumber
% \Xi(\mathbf{x},\bm \varsigma) = \left(x_2-x_1^2\right) \varsigma_1+\left(x_3-x_2^2\right) \varsigma_2-\frac{\varsigma_1^2}{400}-\frac{\varsigma_2^2}{400}+ x_1^2+ x_2^2-2 x_1-2 x_2 + 2 \nonumber
\end{eqnarray}
\indent We select $\bm \varsigma_0 = -1$ to make sure that $G(\bm \varsigma_0) \succeq 0$. Using the Moore-Penrose pseudoinverse, we can obtain the corresponding $\mathbf{x}_0 = G(\bm \varsigma_0)^{-1} F(\bm \varsigma_0)$ $=(0.5,0)$. Taking $(\mathbf{x}_0, \bm \varsigma_0)$ as initial point for the nonlinear equations in (\ref{eq:16}), and after 4 iterations with 0.304920 seconds, we obtain $(\mathbf{x}^{*}, \bm \varsigma^{*}) = (1,1,0)$ and $P(\mathbf{x}^\ast) = 0$.
\\
\indent \textit{Strategy 2}\\
\indent The canonical dual function is
\begin{eqnarray}
P^d(\varsigma) = 1 -\frac{\varsigma^2}{400} - \frac{1}{2} \Big(2,  -\varsigma \Big)
\begin{pmatrix}
2 -2\varsigma & 0\\
0 &0
\end{pmatrix}^{+}
\begin{pmatrix}
2\\
-\varsigma
\end{pmatrix}.
\nonumber
\end{eqnarray}
\indent Due to the singularity of matrix $G(\varsigma)$, we can not get a proper form of $P^d(\varsigma)$; thus it becomes difficult to solve the nonlinear equations in (\ref{eq:17}).\\
\indent \textit{Strategy 3}\\
\indent The same situation happens as above, we can choose some possible initial point $\varsigma_0$ to guarantee $G(\varsigma_0) \succeq 0$, but the calculation of the singular matrix is quite complicated.\\
\indent \textit{Strategy 4}\\
\indent Instead, we choose the same $\varsigma_0$, and using the Moore-Penrose pseudoinverse, we can obtain the corresponding $\mathbf{x}_0 = (0.5,0)$ for $P(\mathbf{x})$.
Taking $\mathbf{x}_0$ as initial point for $P(\mathbf{x})$, we can finally reach $\mathbf{x}^{\ast} = (1,1)$ within 20 iterations and 0.352472 seconds, and then the corresponding $P(1,1) = 0$.\\
\indent Furthermore, we continue to consider the Rosenbrock function in terms of large dimensions. We choose
$\bm \varsigma_0 = (-1,\cdots,-1)$ to guarantee $G(\bm\varsigma_0) \succeq 0$, and then get the corresponding initial point
$\mathbf{x}_0 = (0.5,\underbrace{0.75,\cdots,0.75}_{n-2},0)$ for $P(\mathbf{x})$. General results of the Rosenbrock function by Strategy 4 are given in Table \ref{tab:Rosenbrock}.\\
\begin{table}[!htbp]
\caption{Results of the Rosenbrock function using Strategy 4}
\label{tab:Rosenbrock}       % Give a unique label
\begin{tabular}{ccccc}
\hline\noalign{\smallskip}
n & $\mathbf{x}^{\ast}$ & $P(\mathbf{x}^{\ast})$ & iterations & time(s)  \\
\noalign{\smallskip}\hline\noalign{\smallskip}
2 & (1,1) & 2.0269e-011 & 20 & 0.352472 \\
5 & (1,$\cdots$,1) & 5.4958e-011 & 29 & 0.405747 \\
10 & (1,$\cdots$,1) & 1.0633e-010 & 31 & 0.409724 \\
20 & (1,$\cdots$,1) & 5.3688e-011 & 37 & 0.423663 \\
50 & (1,$\cdots$,1) & 1.6986e-009 & 42 & 0.554678 \\
100 & (1,$\cdots$,1) & 3.7337e-010 & 50 & 0.727062 \\
200 & (1,$\cdots$,1) & 1.5632e-010 & 55 & 1.329283 \\
500 & (1,$\cdots$,1) & 3.0872e-010 & 54 & 3.508815 \\
1000 & (1,$\cdots$,1) & 5.0893e-010 & 56 & 8.763668 \\
2000 & (1,$\cdots$,1) & 3.7200e-010 & 60 & 28.264277 \\
3000 & (1,$\cdots$,1) & 7.3433e-010 & 62 & 57.669020 \\
4000 & (1,$\cdots$,1) & 1.0350e-009 & 61 & 92.344600 \\
5000 & (1,$\cdots$,1) & 1.0340e-009 & 66 & 144.069188\\
\noalign{\smallskip}\hline
\end{tabular}
\end{table}
\\
\indent Compared the results for the Rosenbrock function with those gained by most popular stochastic methods, like PSO (CLPSO, APSO)\citep{liang,zhan} and DE (SaDE)\citep{qin,das}, we can conclude that the strategy used in this paper by canonical duality theory is much more superior. To the best of our knowledge, it is the first time to solve the Rosenbrock function optimization problem up to 5000 dimension in such a short time.\\
\\
\indent \textit{Example 5} (Dixon and Price function)
\begin{eqnarray}
f_5 = (x_1 - 1)^2 + \sum_{i=2}^{n} i (2x_i^2 - x_{i-1})^2. \nonumber
\end{eqnarray}
\indent The landscape of Dixon and Price function is given in Fig. \ref{fig:dp}.
\\
\indent We firstly rewrite it to the standard form, and then we can get $a_k = 2(k+1)$, $A_k = 4 I_{k+1}$, $b_k = -e_k$, $c_k = 0$, $Q = Diag\{2,\underbrace{0,\cdots,0}_{n-1}\}$,
$f = [2,\underbrace{0,\cdots,0}_{n-1}]^T$, where $k = 1,2,\cdots,n-1$, and then $G(\bm \varsigma) = Q + \sum_{k=1}^{n-1} \varsigma_k A_k = Diag\{2,4\varsigma_1,\cdots, 4\varsigma_{n-1}\}$, $F(\bm \varsigma) = f - \sum_{k=1}^{n-1} \varsigma_k b_k = [2+\varsigma_1,\varsigma_2,\cdots,\varsigma_{n-1},0]$.\\
\indent It is not difficult to find that, for any dual feasible solution $\bm\varsigma^{*}$, if we substitute back then we will find that the last component of the corresponding $\mathbf{x}^{*}$
will always be zero, which indicates that the first three strategies will be invalid. However, the fourth strategy can still survive if we make some minor revisions. For simplicity, we choose $\bm\varsigma_0 = (2,\cdots,2)$ to make sure $G(\bm\varsigma_0) \succeq 0$, and then get the corresponding initial point
$\mathbf{x}_0 = (2,\underbrace{0.25,\cdots,0.25}_{n-2},0)$ for $P(\mathbf{x})$. We don't use the $\mathbf{x}_0$ directly but translate the point to $\mathbf{x}_0 = \mathbf{x}_0 + 1 = (3,\underbrace{1.25,\cdots,1.25}_{n-2},1)$. Taking the revised initial point for the primal problem, the general results of the Dixon and Price function are given in Table \ref{tab:DixonPrice}.
\begin{table}[!htbp]
\caption{Results of the Dixon and Price function using Strategy 4}
\label{tab:DixonPrice}       % Give a unique label
\begin{tabular}{ccccc}
\hline\noalign{\smallskip}
n & $\mathbf{x}^{\ast}$ & $P(\mathbf{x}^{\ast})$ & iterations & time(s)  \\
\noalign{\smallskip}\hline\noalign{\smallskip}
2 & (1,0.7071) & 3.1388e-015 & 12 & 0.213785 \\
5 & (1,$\cdots$,0.5221) & 8.4890e-014 & 21 &  0.206739 \\
10 & (1,$\cdots$,0.5007) & 5.4620e-012 & 30 & 0.218370 \\
20 & (1,$\cdots$,0.5000) & 9.1666e-011 & 46 & 0.245217 \\
50 & (1,$\cdots$,0.5000) & 3.4299e-010 & 79 & 0.388959 \\
100 & (1,$\cdots$,0.5000) & 3.6424e-009 & 108 & 0.757873 \\
200 & (1,$\cdots$,0.5000) & 1.0303e-008 & 154 & 1.720907 \\
500 & (1,$\cdots$,0.5000) & 3.1588e-008 & 242 & 7.814894 \\
1000 & (1,$\cdots$,0.5000) & 6.8696e-008 & 342 & 28.862242 \\
2000 & (1,$\cdots$,0.5000) & 1.3657e-007 & 480 & 124.977932 \\
3000 & (1,$\cdots$,0.5000) & 2.4159e-007 & 581 & 270.350883 \\
4000 & (1,$\cdots$,0.5000) & 2.2758e-007 & 675 & 526.158263 \\
5000 & (1,$\cdots$,0.5000) & 3.5225e-007 & 747 & 854.212220\\
\noalign{\smallskip}\hline
\end{tabular}
\end{table}
\section{Conclusion}
To efficiently apply the canonical duality theory for real world problems, four strategies
are proposed to develop algorithms based on the theory. The former two strategies should calculate the
staionary points, in other words, solving nonlinear equations, while the later strategies use numerical optimization algorithms based on unconstrained methods. Some experimental results are given to illustrate the details of using the four strategies for fourth-order polynomial benchmark functions,
and we find that various strategies have different degrees of complexity. To some extent, the canonical duality theory can eliminate the gap between deterministic and
stochastic methods. In our future work, we will try to use stochastic methods to design efficient algorithms for the powerful canonical duality theory.
\section*{Acknowledgments}
Xiaojun Zhou's research is supported by China Scholarship Council, and Chunhua Yang is
supported by the National Science Found
for Distinguished Young Scholars of China (Grant No.
61025015).

\end{document}